\documentclass[oneside]{amsart}
\usepackage[utf8]{inputenc}
\usepackage{amsmath,amstext,amsopn,amssymb, amsfonts, amsthm, mathtools, float}
\usepackage[dvipsnames]{xcolor}
\usepackage{geometry}
\usepackage{bbm}
\usepackage{pdfsync}
\usepackage{tikz}
\usepackage{enumitem}
\usepackage{subfig}
\usepackage{caption}
\usepackage{graphicx}
\usepackage{mathrsfs}
\usepackage[ruled,vlined,algosection]{algorithm2e}
\DeclareMathAlphabet{\mathpzc}{OT1}{pzc}{m}{it}
\usepackage[pdftex,bookmarks,colorlinks,breaklinks]{hyperref}  % PDF hyperlinks, with coloured links
\definecolor{dullmagenta}{rgb}{0.4,0,0.4}   % #660066
\definecolor{darkblue}{rgb}{0,0,0.4}
\definecolor{darkgreen}{rgb}{0,0.4,0}
\hypersetup{linkcolor=darkblue,citecolor=blue,filecolor=dullmagenta,urlcolor=darkblue} % coloured links

{}

\newtheorem{definition}{Definition}
\newtheorem*{definition*}{Definition}

\newtheorem{theorem}{Theorem}
\newtheorem*{theorem*}{Theorem}

\newtheorem*{conjecture*}{Conjecture}

\newtheorem{question}{Question}
\newtheorem*{question*}{Question}

\newtheorem*{lemma*}{Lemma}

\newtheorem{corollary}[theorem]{Corollary}
\newtheorem*{corollary*}{Corollary}

\newtheorem*{remark*}{Remark}
\numberwithin{equation}{section}
\numberwithin{theorem}{section}

\makeatletter
\newcommand{\customlabel}[2]{%
   \protected@write \@auxout {}{\string \newlabel {#1}{{#2}{\thepage}{#2}{#1}{}} }%
   \hypertarget{#1}{#2}
}
\makeatother

\def\Xint#1{\mathchoice%
   {\XXint\displaystyle\textstyle{#1}}%
   {\XXint\textstyle\scriptstyle{#1}}%
   {\XXint\scriptstyle\scriptscriptstyle{#1}}%
   {\XXint\scriptscriptstyle\scriptscriptstyle{#1}}%
   \!\int}
\def\XXint#1#2#3{{\setbox0=\hbox{$#1{#2#3}{\int}$}
     \vcenter{\hbox{$#2#3$}}\kern-.5\wd0}}

\def\dashint{\Xint-}

\newcommand{\sh}{\operatorname{sh}}

\usepackage{color}

\hypersetup{colorlinks = true, urlcolor = blue, linkcolor = blue, citecolor = blue}

\usepackage{setspace}
\date{\today}

\makeindex
\begin{document}

\title[weak-type maximal function estimates on the infinite-dimensional torus]{weak-type maximal function estimates \\ on the infinite-dimensional torus}

\author{Dariusz Kosz}
	\address{Dariusz Kosz (\textnormal{dkosz@bcamath.org}) \newline
	Basque Center for Applied Mathematics, 48009 Bilbao, Spain \newline 
	Wroc{\l}aw University of Science and Technology, 50-370 Wroc{\l}aw, Poland		
	}

\author{Guillermo Rey}
	\address{Guillermo Rey (\textnormal{guillermo.rey@uam.es}) \newline
	Universidad Aut{\'o}noma de Madrid, 28049 Madrid, Spain		
	}

\author{Luz Roncal}	
\address{Luz Roncal (\textnormal{lroncal@bcamath.org})\newline
Basque Center for Applied Mathematics,
48009 Bilbao, Spain\newline Ikerbasque, Basque Foundation for Science, 48011 Bilbao, Spain\newline Department of Mathematics\\
UPV/EHU\\
Apto. 644, 48080 Bilbao, Spain}

	\begin{abstract} We prove necessary and sufficient conditions for the weak-$L^p$ boundedness, for $p \in (1,\infty)$, of a maximal operator on the infinite-dimensional torus.
    In the endpoint case $p=1$ we obtain the same weak-type inequality enjoyed by the strong maximal function in dimension two. Our results are quantitatively sharp. 
		
		
	\end{abstract}

\maketitle

\section{Introduction} \label{S1}

The infinite-dimensional torus $\mathbb{T}^\omega$ is the product of countably many copies of the one-dimensional torus $\mathbb{T}$.
Motivated by the interest in the study of Fourier analysis on $\mathbb{T}^\omega$ pointed out by Rubio de Francia \cite{RdF78, RdF80},
in recent times there has been an effort to understand how the fundamental tools of classical harmonic analysis look like in the setting of $\mathbb{T}^\omega$ \cite{FR20, Ko21, Ko22, KMPRR22}.
As noticed in \cite[Theorem~1.1]{Ko22}, by the results of \cite{CW71}, \cite{PS09}, and \cite{Ed98}, topologically infinite-dimensional spaces cannot satisfy the doubling condition,
which rules out the application of the general theory of spaces of homogeneous type in the Coifman--Weiss \cite{CW71} sense to the objects defined on $\mathbb{T}^{\omega}$. 
What this observation means is that new tools and ideas are needed to develop a theory of harmonic analysis in $\mathbb{T}^\omega$. 

Recent research has focused on obtaining an analogue of the  Calder\'on--Zygmund decomposition and the introduction of a suitable ``dyadic'' basis $\mathcal{R}_0$
and a ``dyadic'' Hardy--Littlewood maximal function $\mathcal M_{\mathcal{R}_0}$ \cite{FR20} which satisfies the weak-type $(1,1)$ inequality.
In \cite{Ko21} it was shown that $\mathcal{M}_{\mathcal{R}}$ does not satisfy any restricted weak-type $(p,p)$ estimate, for $p \in [1, \infty)$, where $\mathcal{R}$ is the wider natural basis, also introduced in \cite{FR20}.
The intermediate case was studied in \cite{KMPRR22}, where bases $\mathcal{R}_0\subset\mathcal{S}\subset \mathcal{R}$
were found for which the corresponding maximal functions present a satisfactory control on Lebesgue spaces. 

In particular, in \cite{KMPRR22} the authors managed to construct explicit bases $\mathcal{S}$ for which the associated maximal function $\mathcal{M}_{\mathcal{S}}$
has an arbitrary restricted weak-type behavior. The relevant result can be summarized as follows:

\begin{theorem}[{\cite[Corollary~1.5]{KMPRR22}}]\label{thm:KMPRR22}
Let $p_0 \in [1, \infty)$.
\begin{enumerate}
    \item It is possible to find a basis $\mathcal{S}$ with $\mathcal{R}_0\subseteq\mathcal{S}\subset \mathcal{R}$ such that $\mathcal M_{\mathcal{S}}$ is of restricted weak-type $(p,p)$ if $p \in [p_0, \infty]$ but is not if $p \in [1, p_0)$.
    \item Similarly, it is possible to find a basis $\mathcal{S}$ with $\mathcal{R}_0\subset\mathcal{S}\subset \mathcal{R}$ such that $\mathcal M_{\mathcal{S}}$ is of restricted weak-type $(p,p)$ if $p \in (p_0, \infty]$ but is not if $p \in [1, p_0]$.
\end{enumerate} 

\end{theorem}

We recall that Theorem \ref{thm:KMPRR22} was obtained through a special construction involving the so-called  \textit{$(\varepsilon,d)$-configurations around sets $Q\in \mathcal{R}_0$}, which are defined in Subsection \ref{sub:maximal}. The approach in \cite{KMPRR22} made use of a combination of elements from harmonic analysis and probability. 

The main purpose of this note is to improve the results in \cite{KMPRR22} by showing sharp weak-type $(p,p)$ estimates,
for $p\in (1,\infty)$, and a sharp endpoint estimate, for $p=1$, for the maximal operator $\mathcal M_{\mathcal{S}}$,
where $\mathcal{S}$ is a countable union of disjointly placed $(\varepsilon,d)$-configurations around sets. 

%%%%%%%%%%%%%%%%%%%%%%%%%%%%%%%%%%
\subsection{Maximal functions and $(\varepsilon,d)$-configurations}
\label{sub:maximal}
%%%%%%%%%%%%%%%%%%%%%%%%%%%%%%%%%%

Given a basis $\mathcal B$, that is, a collection of measurable sets $B \subset \mathbb T^{\omega}$, we define the maximal operator associated with $\mathcal B$ by
\begin{align*}
\mathcal{M}_{\mathcal{B}}f \coloneqq \sup_{B \in \mathcal{B}} \frac{\mathbbm{1}_B}{|B|}\int_B |f| \, dx,
\end{align*}
where $dx$ is the normalized Haar measure on $\mathbb{T}^{\omega}$, the symbol $\mathbbm{1}_B$ stands for the indicator function of $B$, and $|B| \coloneqq \int \mathbbm{1}_B \, dx$. We always assume that $|B| > 0$ for all $B \in \mathcal B$. Moreover, we use the convention $\mathcal{M}_{\mathcal{B}}f(x) = 0$, whenever there is no $B$ containing $x$. 

Historically, the two most important maximal operators in the context of ${\mathbb T}^{\omega}$ are those associated with the Rubio de Francia bases, classical $\mathcal{R}$ and ``dyadic'' $\mathcal{R}_0$. They are introduced with the aid of specifically chosen subsets of $\mathbb T^{\omega}$ (see \cite{FR20} or \cite{KMPRR22} for precise definitions): 
\begin{itemize}
	\item subgroups $H_1 \subset H_2\subset \cdots$ satisfying $[H_{m+1}:H_m]=2$ and $\mathbb{T}^\omega=\overline{\bigcup_{m\in \mathbb{N}}H_m}$,
	\item open sets $V_1 \supset V_2 \supset \cdots$ corresponding to the quotient groups ${\mathbb T}^{\omega} / H_m$.  
\end{itemize} 
Table~\ref{Tab} shows how the first few objects look like. We use the auxiliary sets $R_n \coloneqq \{0, \frac{1}{n}, \dots, \frac{n-1}{n}\}$, $n \in {\mathbb N}$, and by ${\bf 0}^{(n}$ we mean the zero vector in $\mathbb T^{\omega}$ with the first $n$ coordinates removed.

\begin{table}[H]
	\caption{Objects used to define the Rubio de Francia bases.  }
	\begin{center}
		\begin{tabular}{ r l l }
			\hline
			$m$ & $H_m$ & $V_m$ \\
			\hline
			$1$ & $R_2 \times \{{\bf 0}^{(1}\}$
			& $(0,\frac{1}{2}) \times {\mathbb T}^{1, \omega}$ \\ 
			$2$ & $R_2 \times R_2 \times \{{\bf 0}^{(2}\}$ 
			& $(0,\frac{1}{2}) \times (0,\frac{1}{2}) \times {\mathbb T}^{2, \omega}$ \\ 
			$3$ & $R_4 \times R_2 \times \{{\bf 0}^{(2}\}$ & $(0,\frac{1}{4}) \times (0,\frac{1}{2}) \times {\mathbb T}^{2, \omega}$ \\
			$4$ & $R_4 \times R_4 \times \{{\bf 0}^{(2}\}$ & $(0,\frac{1}{4}) \times (0,\frac{1}{4}) \times {\mathbb T}^{2, \omega}$ \\ 
			$5$ & $R_4 \times R_4 \times R_2 \times \{{\bf 0}^{(3}\}$ & $(0,\frac{1}{4}) \times (0,\frac{1}{4}) \times (0,\frac{1}{2}) \times{\mathbb T}^{3, \omega}$ \\ 
			$6$ & $R_4 \times R_4 \times R_4 \times \{{\bf 0}^{(3}\}$ & $(0,\frac{1}{4}) \times (0,\frac{1}{4}) \times (0,\frac{1}{4}) \times{\mathbb T}^{3, \omega}$ \\
			$7$ & $R_8 \times R_4 \times R_4 \times \{{\bf 0}^{(3}\}$ & $(0,\frac{1}{8}) \times (0,\frac{1}{4}) \times (0,\frac{1}{4}) \times{\mathbb T}^{3, \omega}$ \\ 
			$8$ & $R_8 \times R_8 \times R_4 \times \{{\bf 0}^{(3}\}$ & $(0,\frac{1}{8}) \times (0,\frac{1}{8}) \times (0,\frac{1}{4}) \times{\mathbb T}^{3, \omega}$ \\ 
			$9$ & $R_8 \times R_8 \times R_8 \times \{{\bf 0}^{(3}\}$ & $(0,\frac{1}{8}) \times (0,\frac{1}{8}) \times (0,\frac{1}{8}) \times{\mathbb T}^{3, \omega}$ \\
			$10$ & $R_8 \times R_8 \times R_8 \times R_2 \times \{{\bf 0}^{(4}\}$ & $(0,\frac{1}{8}) \times (0,\frac{1}{8}) \times (0,\frac{1}{8}) \times (0,\frac{1}{2}) \times{\mathbb T}^{4, \omega}$ \\
			$\cdots$ & $\cdots$ & $\cdots$ \\
			\hline
		\end{tabular}
		\label{Tab}
	\end{center}
\end{table}  

\noindent According to the notation above the Rubio de Francia bases are defined as follows
\[
\mathcal{R}_0 \coloneqq \{ t + V_m : m \in {\mathbb N}, \, t \in H_m \} 
\quad {\rm and} \quad
\mathcal{R} \coloneqq \{ t + V_m : m \in {\mathbb N}, \, t \in {\mathbb T}^\omega \}.
\]

The behavior of $\mathcal M_{\mathcal R_0}$ and $\mathcal M_{\mathcal R}$ acting on $L^p$ functions is well known.
Precisely, $\mathcal M_{\mathcal R_0}$ is not of strong-type $(1,1)$ but it is of weak-type $(1,1)$, as was proven in \cite{FR20} using the classical martingale-type argument,
hence it is of strong-type $(p,p)$ for all $p \in (1, \infty]$.
On the other hand, $\mathcal M_{\mathcal R}$ does not satisfy any strong-type (or weak-type) $(p,p)$ inequality for any finite $p$, as was shown in \cite{Ko21}. 

The observations above prompted the authors of \cite{KMPRR22} to study intermediate bases $\mathcal{R}_0\subset\mathcal{S}\subset \mathcal{R}$ and examine mapping properties of the associated maximal operators. As a result of these investigations, the $(\varepsilon, d)$-configurations around sets were introduced.
 
\begin{definition}\label{epsn}
Let $\varepsilon \in (0, \frac{1}{2}]$ and $d \in \mathbb{N}$. 
 Let $Q \in \mathcal{R}_0$ be such that $Q = V_m + t$ for some $t \in H_m$ with $m \in \mathbb N$ so large that $V_{m}$ has at least $d$ nonfree coordinates. For each $k \in \{1, \dots, d\}$ set $Q_{k} = T_{k} + Q$, where $T_{k} = (0, \dots, 0, (1-\varepsilon)|\pi_k(Q)|, 0, \dots)$, with $\pi_k$ being the projection onto the $k$-th coordinate. Then $\{Q_{1}, \dots, Q_{d}\}$ is called an $(\varepsilon, d)$-configuration around $Q$.
\end{definition}

Observe that a collection $\{Q_{1}, \dots, Q_{d}\}$ which is an $(\varepsilon, d)$-configuration around $Q \in \mathcal R_0$ satisfies the following properties:
\begin{enumerate}
	\item $Q_k \in \mathcal R$ for each $k$,
    \item $|Q_{k}|=|Q|$ and $| Q_{k} \cap Q| = \varepsilon |Q|$ hold for each $k$,
    \item the sets $Q_{k} \setminus Q$, $k \in \{1, \dots, d\}$, are disjoint and of sizes comparable to that of $Q$,
    \item the sets $Q_{k} \cap Q$, $k \in \{1, \dots, d\}$, resembles a collection of independent events. 
\end{enumerate}

\begin{figure}[H]
\begin{tikzpicture}[scale=0.5]
  \usetikzlibrary {arrows.meta}
  \filldraw[White, draw=black, opacity=0.4] (0,0) rectangle (4,4);
  
  \filldraw[NavyBlue, draw=black, opacity=0.4] (3,0) rectangle (7,4);
  \filldraw[NavyBlue, draw=black, opacity=0.4] (0,3) rectangle (4,7);
  \draw[|<->|] (3,-0.5) to (4,-0.5);
  \node at (3.5,-1.25) {\small $\varepsilon |\pi_1(Q)|$};
  \node at (1,1) {\small $Q$};
  \node at (7.5,3) {\small $Q_1$};
  \node at (3,7.5) {\small $Q_2$};
\end{tikzpicture}
  \label{fig:M1}
  \caption{An $(\varepsilon,d)$-configuration with $d=2$. The blue squares are the two elements $Q_1$ and $Q_2$ in the configuration, and the white square at the center is $Q$.}
\end{figure}
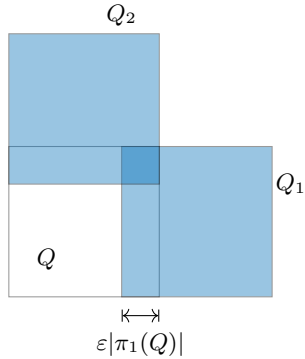

The class of bases considered in \cite{KMPRR22} consists of all collections $\mathcal{S}$ of the form
\begin{align*}
  \mathcal{S} = \bigcup_j \mathcal{S}_j,
\end{align*}
where $\mathcal{S}_j = \{ Q_{j,i}\}_{i=1}^{d_j}$ are $(\varepsilon_j, d_j)$-configurations around some rectangles $Q_{j,0} \subset \mathbb{T}^\omega$ with parameters $\varepsilon_j \in (0, \frac{1}{2}]$ and $d_j \in \mathbb{N}$, and such that the sets $\sh(\mathcal{S}_j)$ are disjoint. Here, for any collection of sets $\mathcal{E}$, 
\begin{align*}
  \sh(\mathcal{E}) \coloneqq \bigcup_{E \in \mathcal{E}} E.
\end{align*}
We shall say that the configurations $\mathcal{S}_j$ are disjoint if the sets $\sh(\mathcal{S}_j)$ are disjoint. 

We noted before that the sets $Q_{j,1}, \dots, Q_{j,d_j}$, when intersected with $Q_{j,0}$, represent independent events.
Indeed, define the probability measure $\mu_j = \frac{\mathbbm{1}_{Q_{j,0}}}{|Q_{j,0}|}$.
Then, for each $\mathcal{E} \subset \mathcal{S}_j$ with $m$ elements,
\begin{align*}
  \mu_j \Big( \bigcap_{Q \in \mathcal{E}} Q \Big) = \varepsilon_j^{m} = \prod_{Q \in \mathcal{E}} \mu_j (Q). 
\end{align*}

In this article, instead of studying restricted weak-type inequalities,
we shall be interested in finding precise bounds for the operator norm of $\mathcal{M}_{\mathcal{S}}$ acting from $L^p$ to $L^{p, \infty}$.
We also deal with the endpoint case $p=1$, which turns out to be quite interesting.
Here we obtain a certain weak-type Orlicz bound similar to the weak-type bound
for the strong maximal function near $L^1$.

To begin with, let us notice that for disjoint configurations $\mathcal{S}_j$ we have
\begin{align*}
  \|\mathcal{M}_{\mathcal{S}}\|_{L^p \to L^{p,\infty}} = \sup_j \|\mathcal{M}_{\mathcal{S}_j}\|_{L^p \to L^{p,\infty}},
\end{align*}
and similarly for other norms.
This allows us to deal with these configurations one at a time.

Let $\mathcal{S}_0$ be a single $(\varepsilon, d)$-configuration around a rectangle $Q_0$, and let $\{Q_i\}_{i=1}^d$ be the translated rectangles in the configuration.
Testing with $f = \mathbbm{1}_{Q_0}$ we obtain
\begin{align*}
  |\{x \in \mathbb{T}^\omega:\, \mathcal{M}_{\mathcal{S}_0}f(x) \geq \varepsilon\}| = |\sh(\mathcal{S}_0)| \geq \frac{d}{2} |Q_0|,
\end{align*}
which gives
\begin{align*}
  \|\mathcal{M}_{\mathcal{S}_0}\|_{L^p \to L^{p,\infty}} \gtrsim \varepsilon d^{\frac{1}{p}}.
\end{align*}
Hence, if $\mathcal{S}$ is a union of disjoint $(\varepsilon_j, d_j)$-configurations, then
\begin{align*}
  \|\mathcal{M}_{\mathcal{S}}\|_{L^p \to L^{p,\infty}} \gtrsim \sup_j \varepsilon_j d_j^{\frac{1}{p}}.
\end{align*}
Therefore, this quantity determines a necessary condition for the weak-$L^p$ boundedness of $\mathcal{M}_{\mathcal{S}}$.
We will denote $A_p(\mathcal{S}_0) \coloneqq \varepsilon d^{\frac{1}{p}}$ for a single $(\varepsilon,d)$-configuration $\mathcal{S}_0$,
and $A_p(\mathcal{S}) \coloneqq \sup_j A_p(\mathcal S_j)$ for $\mathcal{S}$ being a union of disjoint $(\varepsilon_j,d_j)$-configurations $\mathcal{S}_j$.

Our first result shows exactly how the norm $\|\mathcal{M}_{\mathcal{S}_0}\|_{L^p \to L^{p,\infty}}$
depends on $p$ and $A_p(\mathcal{S}_0)$.

\begin{theorem} \label{WeakTypeTheoremP>1}
  Let $\mathcal{S}_0$ be a single $(\varepsilon, d)$-configuration. Then, for $p \in (1, \infty)$ and $q \coloneqq \frac{p}{p-1}$,
  \begin{enumerate}
    \item $\|\mathcal{M}_{\mathcal{S}_0}\|_{L^p \to L^{p,\infty}} \lesssim 1$ if $A_p(\mathcal S_0) \leq q e^{-q}$,
    \item $\|\mathcal{M}_{\mathcal{S}_0}\|_{L^p \to L^{p,\infty}} \lesssim \frac{q}{\log(q \slash A_p(\mathcal S_0))}$  if $qe^{-q} \leq A_p(\mathcal S_0) \leq q e^{-1}$,
    \item $\|\mathcal{M}_{\mathcal{S}_0}\|_{L^p \to L^{p,\infty}} \lesssim e A_p(\mathcal S_0)$ if $A_p(\mathcal S_0) \geq q e^{-1}$.
  \end{enumerate}
  Moreover, the upper bounds are sharp, up to constants independent of $\varepsilon, d, p$.
\end{theorem}

\noindent In the last case, for technical reasons we put $e A_p(\mathcal S_0)$ instead of $A_p(\mathcal S_0)$. With this choice, for each fixed $p$ the postulated upper bound becomes a continuous increasing function of $A_p(\mathcal S_0)$. 

In view of Theorem~\ref{WeakTypeTheoremP>1}, the condition $A_p(\mathcal{S}) < \infty$, which is necessary for the $L^p \to L^{p, \infty}$ boundedness of $\mathcal M_{\mathcal S}$, is in fact sufficient when $p \in (1, \infty)$.

\begin{corollary}
    Let $\mathcal{S}$ be a union of disjoint  $(\varepsilon_j, d_j)$-configurations $\mathcal{S}_j$. Then, for each $p \in (1, \infty)$, 
    \begin{align*}
        \Big( \|\mathcal{M}_{\mathcal{S}}\|_{L^p \to L^{p,\infty}} < \infty \Big) \iff 
        \Big( A_p(\mathcal{S}) < \infty \Big).
    \end{align*}
In particular, for each given $p_0 \in (1, \infty)$, taking $\varepsilon_j = \frac{1}{j+1}$ and $d_j = \lfloor j^{p_0} \rfloor$ (respectively, $d_j = \lfloor \log(j+2) j^{q_0} \rfloor$), we obtain that the range of boundedness is precisely $[p_0, \infty]$ (respectively, $(p_0, \infty]$).
\end{corollary}

\noindent Moreover, we point out an interesting special case where a logarithmic-type improvement occurs.  
\begin{corollary}
    Let $\mathcal{S}$ be a union of disjoint $(\varepsilon_j, d_j)$-configurations $\mathcal{S}_j$. Assume that for some $p_0 \in (1, 2]$ we have $A_{p_0}(\mathcal{S}) = 1$. Then, for $q_0 \coloneqq \frac{p_0}{p_0-1}$,
    \begin{align*}
        \|\mathcal{M}_{\mathcal{S}}\|_{L^{p_0} \to L^{{p_0},\infty}} \lesssim \frac{q_0}{\log(q_0)}.
    \end{align*}
\end{corollary}

When $p=1$, as we mentioned, the situation deteriorates and we do \emph{not} have weak-$L^1$ boundedness when $A_1(\mathcal{S})$ is finite. Indeed, the condition which is sufficient and necessary is $\tilde A_1(\mathcal S) \coloneqq \sup_{j} d_j < \infty$, see \cite[Theorem~1.4]{KMPRR22}. Precisely, restricted weak-type $(p,p)$ inequalities are characterized there, with $\tilde A_1(\mathcal S) < \infty$ corresponding to $p=1$. In fact, for $p=1$ and $\mathcal{M}_{\mathcal S}$, it can be checked that being of restricted weak-type $(1,1)$ is equivalent to being of weak-type $(1,1)$. More generally, let us mention that the equivalence holds for maximal operators on $\mathbb R^d$ or $\mathbb T^d$ associated with a sequence of integrable convolution kernels (see \cite{Mo05} and \cite[Chapter 4.5, Theorem 5.6]{BS88}) but it is not true for arbitrary sublinear operators. Indeed, in \cite{HJ05} the authors provide an example of a sublinear translation-invariant operator acting on $L^1(\mathbb T)$ which is not of weak-type $(1,1)$, although it is of restricted weak-type $(1,1)$ (that is, of weak-type $(1,1)$ for indicator functions).

%At this point one could ask for weak-$L^1$ boundedness limited to characteristic functions. But then the situation does not change and $\tilde A_1(\mathcal S) < \infty$ is sufficient and necessary. This is true because in the unboundedness part of \cite[Theorem~1.4]{KMPRR22} characteristic functions were used.  

Now, if we instead ask for an $L\log L \to L^{1,\infty}$ result, then the condition $A_1(\mathcal S) < \infty$ is sufficient. In fact, we can even weaken it slightly.

\begin{theorem} \label{endpoint-theorem}
  Let $\mathcal{S}_0$ be a single $(\varepsilon, d)$-configuration. Then, for each $\lambda \in (0, \infty)$, 
  \begin{align*}
    |\mathcal{M}_{\mathcal{S}_0}f > \lambda| \lesssim (1+A(\mathcal S_0)) \int \frac{|f|}{\lambda} \log\Big( e + \frac{|f|}{\lambda}\Big),
  \end{align*}
  where $A(\mathcal S_0) \coloneqq \frac{\varepsilon d}{\log(\varepsilon^{-1})}$ is a weakened version of $A_1(\mathcal S_0)$.
  Moreover, the upper bound is sharp, up to a constant independent of $\varepsilon, d, \lambda$.
\end{theorem}

As before, denoting $A(\mathcal S) \coloneqq \sup_j \frac{\varepsilon_j d_j}{\log(\varepsilon_j^{-1})}$ we see that the condition $A(\mathcal S) < \infty$ is sufficient and necessary for the $L\log L \to L^{1,\infty}$ boundedness of $\mathcal M_{\mathcal S}$.  

\begin{corollary}
Let $\mathcal{S}$ be a union of disjoint  $(\varepsilon_j, d_j)$-configurations $\mathcal{S}_j$. 
Then
\begin{align*}
        \Big( \|\mathcal{M}_{\mathcal{S}}\|_{L \log L \to L^{1,\infty}} < \infty \Big) \iff 
        \Big( A(\mathcal{S}) < \infty \Big).
    \end{align*} 
\end{corollary}

The rest of the paper is organized as follows. In Sections~\ref{S2}~and~\ref{S3} we prove Theorems~\ref{endpoint-theorem}~and~\ref{WeakTypeTheoremP>1}, respectively. We decided to reverse the order of proofs, so that we can begin with the simpler one. In Section~\ref{S4} we collect some final comments and open questions. 

\section{The $L^1$ endpoint} \label{S2}

%We have seen that the constant $A_p = \varepsilon d^{\frac{1}{p}}$ controls the weak-type $(p,p)$ norm of $\mathcal{M}_{\mathcal{S}}$ for $p \in (1, \infty)$. However, the natural endpoint $A_1 = \varepsilon d$ does not control the weak-type $(1,1)$ norm, as the latter quantity is just comparable to $d$. To obtain a more subtle endpoint result we add an extra logarithmic term and, as in the case of the strong maximal function, study the modular inequality 
%\begin{align*}
%  |\mathcal{M}_{\mathcal{S}}f > \lambda| \lesssim \int \frac{|f|}{\lambda} \log\Big( e + \frac{|f|}{\lambda}\Big)
%\end{align*}
%where $\mathcal{S}$ is a single $(\varepsilon, d)$-configuration.

Here we show Theorem~\ref{endpoint-theorem}. The proof will be similar to the proof of the strong maximal theorem in \cite{CF75}, but with a slight probabilistic bent. We abbreviate $\mathcal S_0$ and $A(\mathcal S_0)$ to $\mathcal S$ and $A$.
\begin{proof}
	In order to prove that the best constant must be at least comparable to $1$, take $f = \mathbbm{1}_{Q_1}$ and $\lambda = \frac{1}{2}$. Then $| \{ \mathcal M_{\mathcal S} f > \lambda \} | = |Q_1|$ and $\int \frac{|f|}{\lambda} \log\big( e + \frac{|f|}{\lambda}\big) \sim |Q_1|$, as desired. Similarly, take $f = \mathbbm{1}_{Q_0}$ and $\lambda = \frac{\varepsilon}{2 }$. Then $| \{ \mathcal M_{\mathcal S} f > \lambda \} | \sim d|Q_0|$ and $\int \frac{|f|}{\lambda} \log\big( e + \frac{|f|}{\lambda}\big) \sim \varepsilon^{-1} \log(\varepsilon^{-1}) |Q_0|$.
	  
  Let us now prove that the postulated estimate holds. Take nonnegative $f$ and $\lambda = 1$. We let
  \begin{align*}
    \mathcal{E} = \{Q \in \mathcal{S} :  \langle f \rangle_Q > 1\}.
  \end{align*}
  Observe that without loss of generality we can assume that $\mathcal{E} = \mathcal{S}$ so that
  \begin{align*}
    |\sh(\mathcal{E})| \sim d |Q_0|.
  \end{align*}
  Indeed, if $\mathcal{E} \neq \mathcal{S}$ then $\mathcal{E} = \mathcal{S}'$ for some $(\varepsilon, d')$-configuration $\mathcal S' \subset \mathcal S$ with $d' < d$.
  Thus, the upper bound we are looking for is the maximum of the upper bounds related to all subconfigurations $\mathcal S' \subset \mathcal S$ and the upper bound one can get analyzing the case $\mathcal{E} = \mathcal{S}$. After all, it turns out that the last quantity dominates the other ones because the upper bound from the thesis increases with $A$, and thus it is increasing in $d$ if $\varepsilon$ is fixed.
  
  As mentioned, assume that $\{\mathcal{M}_{\mathcal{S}}f > 1\} = \sh(\mathcal{S})$. We write $f = f_1 + f_2$, where 
  $f_1 \coloneqq f \mathbbm{1}_{Q_0 \cap \{ f \leq \frac{1}{2\varepsilon} \}}$ and $f_2 \coloneqq f - f_1$. Then for each $Q \in \mathcal S$ we have $\langle f_1 \rangle_Q \leq \frac{1}{2}$, so that $\langle f_2 \rangle_Q > \frac{1}{2}$ and consequently 
  \begin{align*}
    |\sh(\mathcal{S})| \leq \sum_{Q \in \mathcal{S}} |Q|
    \leq \sum_{Q \in \mathcal{S}} 2 \int_Q f_2 
    = 2 \int f_2 h_{\mathcal{S}}.
  \end{align*}
  Here we have used the \emph{height function}, defined for any family $\mathcal{F}$ of rectangles as
\begin{align*}
  h_{\mathcal{F}} \coloneqq \sum_{R \in \mathcal{F}} \mathbbm{1}_R.
\end{align*}

  Let $\delta \in (0, 1)$ be a constant to be chosen later, slightly smaller than $1$. Then, for $\widetilde{f} \coloneqq \delta^{-1} f_2$,
  \begin{align*}
    \int f_2 h_{\mathcal{S}} =  \delta\int \widetilde{f} h_{\mathcal{S}} 
     = \delta\int_{Q_0} \widetilde{f} h_{\mathcal{S}} +  \delta\int_{\mathbb T^{\omega} \setminus Q_0} \widetilde{f} h_{\mathcal{S}}.
  \end{align*}
  Consider the first term. We would like to bound $\widetilde{f}h_{\mathcal{S}}$ by a sum of the form $\varphi(\widetilde{f}) + \psi(h_\mathcal{S})$, whenever $h_{\mathcal S} > 0$. Let $c \in (0, 1]$ be a small constant to be chosen later. Set $\psi(x) \coloneqq \varepsilon c^{-1} \exp(c x)$. Then
  \begin{align*}
    ab \leq \varphi(a) + \psi(b)
  \end{align*}
  for all $a, b \in (0, \infty)$, with $\varphi$ being the Legendre transform of $\psi$ given by
  \begin{align} \label{Legendre}
    \varphi(x) \coloneqq \sup_{t \in (0, \infty)}(xt - \psi(t)).
  \end{align}
    With this $\varphi$ we have
  \begin{align*} 
    2 \delta \int_{Q_0} \widetilde{f} h_\mathcal{S}
    \leq 2 \delta \int_{ Q_0 } \varphi(\widetilde{f}) + 2 \delta \varepsilon c^{-1} \int_{Q_0} \exp(c h_\mathcal{S}). 
  \end{align*}

  Next, we obtain an exponential-integrability estimate for $h_{\mathcal{S}}$. Here we can use the independence of the functions $\mathbbm{1}_{Q \cap Q_0}$ for $Q \in \mathcal{S}$ with respect to the probability measure $d \mu = \frac{\mathbbm{1}_{Q_0}}{|Q_0|}dx$. Precisely,
  \begin{align*}
    \int_{Q_0} \exp(c h_{\mathcal{S}}) \leq |Q_0|\int\exp(c h_{\mathcal{S}}) \,d\mu
    = |Q_0| \prod_{Q \in \mathcal{S}} \int \exp(c \mathbbm{1}_Q) \, d\mu.
  \end{align*}
  Let $\Lambda \coloneqq \exp(c)$. For every $Q \in \mathcal{S}$ we have
  \begin{align*}
    \int \exp(c \mathbbm{1}_Q) \, d\mu = \frac{\Lambda |Q_0 \cap Q| + |Q_0 \setminus Q|}{|Q_0|} = 1 - \varepsilon + \Lambda \varepsilon = 1 + (\Lambda-1)\varepsilon.
  \end{align*}
  From this, using the estimate $1 + t \leq \exp t$ for $t \in (0, \infty)$, we deduce that
  \begin{align*}
    \int_{Q_0} \exp(c h_{\mathcal{S}}) \leq |Q_0| (1+(\Lambda-1)\varepsilon)^d 
   \leq |Q_0| \exp(d(\Lambda-1)\varepsilon).
  \end{align*}
  Choose $c$ so that $(1+A)(\Lambda-1) = (1+A) (\exp(c) - 1)= 1$. Then 
  \begin{align*}
    \frac{1}{|Q_0|} \int_{ Q_0} \exp(c h_{\mathcal{S}}) \leq\exp( d(\Lambda-1) \varepsilon)
    = \exp( d(1+A)^{-1} \varepsilon)
    \leq \varepsilon^{-1}.
  \end{align*}
  Consequently, since $c=\log(1+(1+A)^{-1})$ and $\log(1 + t) \sim t$ for $t \in (0,1)$,
   \begin{align*}
  2 \delta \varepsilon c^{-1} \int_{Q_0} \exp(c h_\mathcal{S})
  \leq 2 \delta c^{-1} |Q_0| \sim \delta \big(1+ \tfrac{\varepsilon d}{ \log(\varepsilon^{-1})} \big) |Q_0| \sim \delta d |Q_0|.
  \end{align*}
  This means that after a suitable choice of $\delta$ comparable to $1$ we have
  \begin{align*}
  2 \delta \varepsilon c^{-1} \int_{Q_0} \exp(c h_\mathcal{S})
  \leq \tfrac{1}{2} |\sh(\mathcal{S})| 
  \end{align*}
  so that
  \begin{align*}
  |\sh(\mathcal{S})| \leq 4 \delta \int_{ \mathbb T^{\omega} \setminus Q_0} \widetilde{f} h_{\mathcal{S}}
  + 4 \delta \int_{ Q_0 } \varphi(\widetilde{f}).
  \end{align*}
  
  The first term above is easy to handle. 
  Indeed, since $h_{\mathcal S} \leq 1$ outside $Q_0$, we have
  \begin{align*}
  4 \delta \int_{ \mathbb T^{\omega} \setminus Q_0} \widetilde{f} h_{\mathcal{S}} \leq 4 \delta  \int \widetilde{f} \log (e + \widetilde{f}).
  \end{align*}
  It remains to focus on $\varphi$. The supremum in \eqref{Legendre} is achieved when $t = \frac{1}{c} \log(\frac{x}{\varepsilon})$. This gives
  \begin{align*}
    \varphi(x) &\leq \tfrac{x}{c} \log(\tfrac{x}{\varepsilon})
    \leq \tfrac{x}{c} \big( \log(e+x) + \log(e + \tfrac{1}{\varepsilon}) \big)
    \lesssim \tfrac{x}{c} \log(e + x), 
    %\lesssim (1+A) x \log(e + x).
  \end{align*}
  provided that $x \gtrsim \frac{1}{\varepsilon}$.
   Since $c \sim (1+A)^{-1}$, $\delta \sim 1$, and $f_2 > \frac{1}{2\varepsilon}$ on $Q_0 \cap \{ f > 0 \}$, we arrive at
  \begin{align*}
    |\sh(\mathcal{S})| &\lesssim (1 + \tfrac{1}{c}) \int \delta^{-1} f\log(e+ \delta^{-1} f)
    \lesssim (1+A) \int f \log(e+f),
  \end{align*}
  as desired.
\end{proof}

\section{weak-type $(p,p)$ estimates for $p \in (1, \infty)$} \label{S3}

Here we show Theorem~\ref{WeakTypeTheoremP>1}. To motivate the argument, consider a nonnegative function $f$ and $\lambda  \in (0, \infty)$. Then, writing shortly $\mathcal S$ instead of $\mathcal S_0$, we obtain
\begin{align*}
  \{\mathcal{M}_{\mathcal{S}}f > \lambda \} = \sh(\mathcal{E}),
\end{align*}
where  $\mathcal{E} \coloneqq \{R \in \mathcal{S} : \int_R f > \lambda |R|\}$.
As before, we can estimate $|\sh(\mathcal{E})|$ using the height function
\begin{align*}
  |\{\mathcal{M}_{\mathcal{S}}f > \lambda \}| = |\sh(\mathcal{E})| \leq \sum_{R \in \mathcal{E}} |R| \leq \frac{1}{\lambda} \sum_{R} \int_R f = \frac{1}{\lambda}\int f h_{\mathcal{E}}.
\end{align*}
%The problem is that the rectangles in $\mathcal{E}$ are not pairwise disjoint, so the sum of characteristic functions of the elements in $\mathcal{E}$ is not uniformly bounded.
%This is the same problem that one encounters when studying weak-type estimates for the strong maximal function at the $L^1$ scale, and indeed we can use some of the techniques developed there.

%For any family $\mathcal{F}$ of rectangles define its \emph{height function} $h_{\mathcal{F}}$ as
%\begin{align*}
%  h_{\mathcal{F}} \coloneqq \sum_{R \in \mathcal{F}} \mathbbm{1}_R.
%\end{align*}
%With this notation, we have
%\begin{align*}
%  |\{\mathcal{M}_{\mathcal{S}}f > \lambda \}| \leq \frac{1}{\lambda} \int f h_{\mathcal{E}}.
%\end{align*}
Then, following the steps in \cite{CF75}, we notice that it is possible to derive weak-type estimates for $\mathcal{M}_{\mathcal{S}}$ from Lebesgue norm estimates for $h_{\mathcal{\mathcal{E}}}$.
For example, if one can prove
\begin{align*}
  %=\|h_{\mathcal{E}}\|_{L^q}^q 
  \int h_{\mathcal{E}}^q 
  \leq C |\sh(\mathcal{E})|,
\end{align*}
for $q = \frac{p}{p-1}$ and some $C \in (0, \infty)$, then by H\"older's inequality
\begin{align*}
  |\{\mathcal{M}_{\mathcal{S}}f > \lambda \}| = |\sh(\mathcal{E})| &\leq \lambda^{-1} \|f\|_{L^p} C^{\frac{1}{q}} |\sh(\mathcal{E})|^{\frac{1}{q}},
\end{align*}
which after rearranging terms becomes
\begin{align} \label{WeakTypeEstimatePre}
  \lambda |\{\mathcal{M}_{\mathcal{S}}f > \lambda \}|^{\frac{1}{p}} \leq C^{\frac{1}{q}} \|f\|_{L^p}.
\end{align}
These ideas are the first ingredient in the proof of Theorem~\ref{WeakTypeTheoremP>1}.
In order to compute the $L^q$ norm of the height function we also need ideas from probability, in particular we will use the following
theorem of Lata{\l }a \cite{La97}.
\begin{theorem}[{\cite[Theorem~1]{La97}}] \label{LatalaThm}
  Let $X_1, \dots, X_n$ be independent nonnegative random variables. Then, for $p \in (1, \infty)$,
  \begin{align*}
    \|X_1 + \dots + X_n\|_{L^p} \sim \|(X_i)\|_{L^p},
  \end{align*}
  where the implicit constant does not depend on $p$ or $n$, and
  \begin{align*}
    % \|(X_i)\|_{L^p} := \inf\Big\{ \lambda > 0:\, \prod_{m=1}^n \int \Big( 1 + \frac{X_i}{\lambda} \Big)^p  \leq e^p \Big\}.
    \|(X_i)\|_{L^p} \coloneqq \inf\Big\{ \eta  \in (0, \infty) : \prod_{i=1}^n \mathbb{E} \Big( 1 + \frac{X_i}{\eta} \Big)^p  \leq e^p \Big\}.
  \end{align*}
\end{theorem}

With these tools, we are ready to prove Theorem~\ref{WeakTypeTheoremP>1}.
\begin{proof}[Proof Theorem~\ref{WeakTypeTheoremP>1}]
  Fix $p \in (1, \infty)$ and let $q = \frac{p}{p-1}$. In order to obtain the desired maximal inequality, we will be interested in estimating $L^q$ norms of $h_{\mathcal{E}}$, where $\mathcal{E}$ is a subcollection of a single $(\varepsilon, d)$-configuration around a rectangle $Q_0$. 
  
  As before, without loss of generality we can assume that $\mathcal{E} = \mathcal{S}$ so that
  \begin{align*}
    |\sh(\mathcal{E})| \sim d |Q_0|.
  \end{align*}
  This follows because the upper bound from the thesis increases with $A_p(\mathcal S)$, and thus it is increasing in $d$ if $\varepsilon$ is fixed.

  First, we observe that the function $h_{\mathcal{E}}$ equals $1$ on $\sh(\mathcal{E}) \setminus Q_0$, which gives the sharp bound
  \begin{align*}
    \int h_{\mathcal{E}}^q
    %\|h_{\mathcal{E}}\|_{L^q}^q 
    = (1-\varepsilon) d |Q_0| + \int_{Q_0} h_{\mathcal{E}}^q  \leq d |Q_0| + |Q_0| \int_{Q_0} h_{\mathcal{E}}^q \,d\mu,
  \end{align*}
  where $\mu$ is the probability measure obtained by restricting $dx$ to $Q_0$, that is, $d \mu \coloneqq |Q_0|^{-1} \mathbbm{1}_{Q_0}dx$.
  Let us now enumerate the intersections of rectangles $Q \in \mathcal E$ with $Q_0$ by
  \begin{align*}
    \{Q \cap Q_0 : Q \in \mathcal{E}\} \eqqcolon \{Q_i\}_{i=1}^d.
  \end{align*}
  Set $f_i := \mathbbm{1}_{Q_i}$. By the discussion above, the functions $\{f_i\}$ are independent with respect to $\mu$, so by Theorem~\ref{LatalaThm} it suffices to
  estimate $\|(f_i)\|_{L^q}$.

  In order to compute $\|(f_i)\|_{L^q}$ it suffices to find the unique parameter $\eta \in (0,\infty)$ such that
  \begin{align*}
    \prod_{i=1}^d \int \Big(1+\frac{f_i}{\eta}\Big)^q = e^q.
  \end{align*}
  Note that
  \begin{align*}
    \prod_{i=1}^d \int \Big(1+\frac{f_i}{\eta}\Big)^q = \Big[ 1 + \varepsilon \Big( \Big(1 + \frac{1}{\eta}\Big)^q-1 \Big)\Big]^d.
  \end{align*}
  Solving this for $\eta$ we arrive at
  \begin{align*}
    \eta = \frac{1}{\big(1 + \frac{\exp(qd^{-1})-1}{\varepsilon}\big)^{1 \slash q} - 1}.
  \end{align*}
  Since $|\sh(\mathcal{E})| \sim d|Q_0|$, this gives \eqref{WeakTypeEstimatePre} with $C^{\frac{1}{q}} \sim 1 + d^{-\frac{1}{q}}\eta$.  Thus, the relevant quantity for us is $d^{-\frac{1}{q}}\eta$ and our goal is to estimate
  the supremum of
  \begin{align} \label{RelevantQuantity}
    \frac{d^{-1 \slash q}}{\big(1 + \frac{\exp(qd^{-1})-1}{\varepsilon}\big)^{1 \slash q} - 1},
  \end{align}
  taken over all pairs $(\varepsilon, d)$ satisfying simultaneously
  \begin{enumerate}
    \item $\varepsilon \in (0, \frac{1}{2}]$, \label{Constraint1}
    \item $\varepsilon d^{\frac{1}{p}} = A$, \label{Constraint2}
    \item $d \in \mathbb{N}$, \label{Constraint3}
  \end{enumerate}
  with $A$ being a fixed positive real number. Using \eqref{Constraint2} and putting $\varepsilon = Ad^{-\frac{1}{p}}$  in \eqref{RelevantQuantity} we arrive at
  \begin{align*}
  \eqref{RelevantQuantity} =
  \frac{d^{-1 \slash q}}{\big(1 + \frac{\exp(qd^{-1})-1}{Ad^{-1/p}}\big)^{1 \slash q} - 1}.
  \end{align*}
  
  Note that \eqref{Constraint1} and \eqref{Constraint2} force $d \geq (2A)^p$, so that the range $d \in \mathbb N \cap [(2A)^p, \infty)$ should be analyzed. However, we will work with $d \in [1,\infty)$ first and only then take this restriction into account.
  
  Suppose that $d \leq q$. Then $e^{-\frac{1}{e}} \leq q^{-\frac{1}{q}} \leq d^{-\frac{1}{q}} \leq 1$ and the numerator in \eqref{RelevantQuantity} is unimportant.
  Furthermore,
  \begin{align*}
    (1-e^{-1}) \exp(qd^{-1}) \leq \exp(qd^{-1})-1 \leq \exp(qd^{-1})
  \end{align*}
  and, for some constant $c_1 \in [1,e^{\frac{1}{e}}]$,
  \begin{align*}
    \frac{\exp(qd^{-1})}{d^{-1 \slash p}} = d\frac{\exp(qd^{-1})}{d^{1 \slash q}} = c_1 d\exp(qd^{-1})
  \end{align*}
  Using these estimates together with $\frac{d}{A} > 0$ and $\frac{1}{q} \in (0,1)$, we conclude that \eqref{RelevantQuantity} equals
  \begin{align} \label{RelevantQuantity1}
  \frac{1}{\big(1 + \frac{d\exp(qd^{-1})}{A}\big)^{1 \slash q} - 1}
  \end{align}
   when $d \in [1,q]$, up to a universal constant independent of $d, q, A$. The function $t \mapsto t \exp(q t^{-1})$, defined for $t \in [1,q]$, has derivative
  \begin{align*}
    \exp(qt^{-1}) + t\exp(qt^{-1}) \frac{-q}{t^2} = \exp(qt^{-1})\Big( 1 - \frac{q}{t} \Big) \leq 0
  \end{align*}
  which means that \eqref{RelevantQuantity1} is increasing in $d$.

  Suppose now that $d \geq q$. Then $\exp(q d^{-1})-1$ is well approximated by $qd^{-1}$, namely,
  \begin{align*}
    qd^{-1} \leq \exp(qd^{-1})-1 \leq (e-1)qd^{-1}.
  \end{align*}
  Thus, for some constant $c_2 \in [1,e-1]$,
  \begin{align*}
    \eqref{RelevantQuantity} = \frac{d^{-1 \slash q}}{\big( 1 + c_2\frac{qd^{-1}}{Ad^{-1 \slash p}}\big)^{1 \slash q}-1}
    = \frac{d^{-1 \slash q}}{\big( 1 + c_2\frac{q}{A}d^{-1 \slash q}\big)^{1 \slash q}-1}.
  \end{align*}
  As before, we conclude that \eqref{RelevantQuantity} equals
  \begin{align} \label{RelevantQuantity2}
  \frac{d^{-1 \slash q}}{\big( 1 + \frac{q}{A}d^{-1 \slash q}\big)^{1 \slash q}-1}
  \end{align}
  when $d \in [q,\infty)$, up to a universal constant independent of $d, q, A$. Also, it is easy to see that \eqref{RelevantQuantity2} is decreasing in $d$ in the considered range.
  
  Now observe that for $d=q$ both \eqref{RelevantQuantity1} and \eqref{RelevantQuantity2} equal
    \begin{align} \label{SimplifiedRQ}
    \frac{1}{\big(1+\frac{q}{A}\big)^{1 \slash q}-1}
    \end{align}
  up to a universal constant, and the same is true for \eqref{RelevantQuantity2} with $d = \lceil q \rceil$. Hence, if $(2A)^p \leq q$, then we are left with the expression \eqref{SimplifiedRQ}. On the other hand, if $(2A)^p > q$, then in order to get the desired bound we should calculate \eqref{RelevantQuantity2} with $d = \lceil (2A)^p \rceil $, or equivalently $d = (2A)^p$, that is,
  \begin{align} \label{SimplifiedRQ2}
  \frac{(2A)^{-p \slash q}}{\big( 1 + \frac{q}{A}(2A)^{-p \slash q}\big)^{1 \slash q}-1}.
  \end{align}

First, let us consider \eqref{SimplifiedRQ} for the whole range $A \in (0, \infty)$. We can distinguish three cases.
  \begin{enumerate}
    \item $A \leq q e^{-q}$. Here we obtain
      \begin{align*}
        \eqref{SimplifiedRQ} \leq \frac{1}{(1+e^q)^{1 \slash q}-1} \lesssim 1.
      \end{align*}
    \item $A \geq q e^{-1}$. Here we use the approximation $(1+x)^{\frac{1}{q}}-1 \sim \frac{x}{q}$, valid for $x \in (0,e]$,
      to obtain
      \begin{align*}
        \eqref{SimplifiedRQ} \sim A.
      \end{align*}
    \item $qe^{-q} \leq A \leq q e^{-1}$. Here we use the approximation 
      $\frac{\log x}{q} \leq x^{\frac{1}{q}}-1 \leq \frac{\log x}{q} \exp\big(\frac{\log x}{q}\big)$,
      valid for $x \in [1, \infty)$.
      When $x = 1 + \frac{q}{A}$ then
      $\log x \leq \log(qA^{-1}) \leq q$
      which gives
      $\frac{\log x}{q} \leq x^{\frac{1}{q}}-1 \leq e \frac{\log x}{q}$,
      so that we obtain
      \begin{align*}
      \eqref{SimplifiedRQ} \sim \frac{q}{\log\big(1+\frac{q}{A}\big)}
      \sim \frac{q}{\log\big(\frac{q}{A}\big)}.
      \end{align*}
  \end{enumerate}

Next, let us analyze what can happen if $(2A)^p > q$. Note that $(2A)^{-\frac{p}{q}} = (2A)^{1-p}$ which gives
\begin{align*}
\eqref{SimplifiedRQ2} =
\frac{(2A)^{1-p}}{\big( 1 + \frac{q}{A}(2A)^{1-p}\big)^{1 \slash q}-1}
\sim 	\frac{(2A)^{1-p}}{\big( 1 + \frac{q}{(2A)^{p}} \big)^{1 \slash q}-1}.
\end{align*}
Then we use the approximation $(1+x)^{\frac{1}{q}}-1 \sim \frac{x}{q}$, valid for $x \in (0,1]$,
to obtain
\begin{align*}
\eqref{SimplifiedRQ2} \sim 
\frac{(2A)^{1-p}}{(2A)^{-p}} \sim A.
\end{align*}
It turns out that this part of the analysis essentially has no impact on the results obtained before. Precisely, in the first case, that is, $A \leq q e^{-q}$, we have $A \lesssim 1$. Thus, we obtain 
\begin{align*}
\|\mathcal{M}_{\mathcal{E}}\|_{L^p \to L^{p,\infty}} \lesssim 1 + d^{-\frac{1}{q}}\eta \lesssim 1,
\end{align*}
regardless of whether we replace $1$ by $A$ or not. In the second case, that is, $A \geq q e^{-1}$, we obtain 
\begin{align*}
\|\mathcal{M}_{\mathcal{E}}\|_{L^p \to L^{p,\infty}} \lesssim 1 + d^{-\frac{1}{q}}\eta \lesssim 1 + A \sim A.
\end{align*}
Finally, in the third case, that is, $q e^{-q} \leq A \leq q e^{-1}$, we have $A \sim \frac{q}{\log (\frac{q}{A})}$, whenever $(2A)^p > q$. Indeed, if the latter holds, then $A > \frac{1}{2} q^{\frac{1}{p}} \geq \frac{1}{2} e^{-\frac{1}{e}} q$. Thanks to this $\frac{q}{A}$ is bounded from both sides by positive universal constants, and hence so is $\frac{q}{A} \frac{1}{\log(\frac{q}{A})}$. Consequently, in this case we obtain
\begin{align*}
\|\mathcal{M}_{\mathcal{E}}\|_{L^p \to L^{p,\infty}} \lesssim 1 + d^{-\frac{1}{q}}\eta \lesssim 1 + \frac{q}{\log\big(\frac{q}{A}\big)} \sim \frac{q}{\log\big(\frac{q}{A}\big)}.
\end{align*}      	

  Let us now verify that for these maximal operators the inequalities obtained above are sharp up to constants independent of $A$ and $q$.
  Take $\mathcal{E}$, a single $(\varepsilon,d)$-configuration, and set
  \begin{align*}
    f = \frac{h_{\mathcal{E}}^{q-1}}{\|h_{\mathcal{E}}\|_{L^{q}}^{q/p}}.
  \end{align*}
  Then $\| f \|_{L^p} = 1$ 
  and the averages $\dashint_Q f$ are equal for all $Q \in \mathcal{E}$. Set $\lambda = \dashint_Q f$ for any $Q \in \mathcal{E}$. Then
  \begin{align*}
    |\mathcal{M}_{\mathcal{E}} f \geq \lambda| = |\sh(\mathcal{E})|.
  \end{align*}
  Also, in this case the Cauchy--Schwarz inequality becomes equality
  \begin{align*}
  \int f h_{\mathcal{E}} = \| f \|_{L^p} \| h_{\mathcal{E}} \|_{L^q} = \| h_{\mathcal{E}} \|_{L^q},
  \end{align*}
so that
  \begin{align*}
    |\sh(\mathcal{E})| \geq \frac{1}{2} \sum_{Q \in \mathcal{E}} |Q| = \frac{1}{2} \sum_{Q \in \mathcal{E}} \lambda^{-1} \int_Q f = \frac{1}{2\lambda} \int f h_{\mathcal{E}} 
    = \frac{1}{2\lambda} \|h_{\mathcal{E}}\|_{L^q}.
  \end{align*}
  Therefore, if we know that
  \begin{align*}
    \int h_{\mathcal{E}}^q \geq \tilde C |\sh(\mathcal{E})|,
  \end{align*}
  for some $\tilde C \in (0, \infty)$, then
  \begin{align*}
  \lambda |\sh(\mathcal{E})|^{\frac{1}{p}} \geq 2^{-1} \tilde C^{\frac{1}{q}}.
  \end{align*}
  This shows that the estimates obtained in the last two cases, that is, $A \geq q e^{-1}$ and $q e^{-q} \leq A \leq q e^{-1}$, are sharp. To obtain the same in the first case, that is, $A \leq q e^{-q}$, one should take $f = \mathbbm{1}_{Q_0}$. 
\end{proof}

\section{Comments and questions} \label{S4}

In \cite{KMPRR22}, it was shown that for $p \in (1, \infty)$ the condition $A_p(\mathcal S) < \infty$ characterizes the restricted weak-type $(p,p)$ boundedness of $\mathcal M_{\mathcal S}$.
Theorem~\ref{WeakTypeTheoremP>1} says that exactly the same condition appears in the context of weak-type $(p,p)$ inequalities. This motivates a very natural question that we do not answer here. 

\begin{question} \label{Q1}
Does the condition $A_p(\mathcal S) < \infty$ characterize the strong-type $(p,p)$ boundedness of $\mathcal M_{\mathcal S}$, that is, the $L^{p} \to L^{p}$ boundedness, when $p \in (1, \infty)$? If not, what is the necessary and sufficient condition for this type of boundedness?
\end{question}

\noindent If $A_p(\mathcal S) < \infty$ is indeed the exact condition, then only the sufficiency part needs to be verified.

The condition $A_1(\mathcal S) < \infty$ looks like a natural endpoint condition regarding some type of boundedness of $\mathcal M_{\mathcal S}$ when $p \to 1$. However, it does not characterize neither $L^1 \to L^{1, \infty}$ nor $L \log L \to L^{1,\infty}$ boundedness, and it seems that the desired domain should lie between $L^1$ and $L \log L$. 

\begin{question} \label{Q2}
Does the condition $A_1(\mathcal S) < \infty$ characterize the $X \to L^{1, \infty}$ boundedness of $\mathcal M_{\mathcal S}$ for some function space $X$? More generally, can this condition be viewed as the endpoint condition regarding the strong-type, weak-type, or restricted weak-type $(p,p)$ boundedness of $\mathcal M_{\mathcal S}$, when $p \to 1$? 
\end{question}

The estimate in Theorem~\ref{endpoint-theorem} resembles the one satisfied by the strong maximal function in $\mathbb{R}^2$.
This is not coincidental, in fact $(\varepsilon,d)$-configurations behave, in a certain sense, like a two-parameter family of rectangles.
Indeed, when intersected with the central rectangle, the intersections have all the same projections (normalized with respect to the projections of the central rectangle) except one,
and thus behave like a family of rectangles in $\mathbb{R}^d$ with exactly one side different than the rest (which side depends on the rectangle). We recall that Zygmund \cite{Zy67} proved that a basis of axis-parallel rectangles with exactly $k$ different sides differentiates $L(\log L)^{k-1}$ regardless of the ambient dimension where these rectangles may lie. Also, in this direction, A. C\'ordoba showed that the maximal operator associated with a basis of axis-parallel rectangles with sidelengths $(s,t,\phi(s,t))$, where $\phi(s,t)$ is monotone in each variable separately, maps not only $L(\log L)^2$ but $L(\log L)$ itself into $L^{1,\infty}$, see e.g. \cite{C79}.

A natural question, which is also a special case of Question~\ref{Q2}, is whether one can replace the logarithm in Theorem~\ref{endpoint-theorem} with a function with slower growth at infinity.
In particular, one can ask for which $\varphi$ do we have an inequality like
\begin{align*}
  |\mathcal{M}_{\mathcal{S}}f > \lambda| \lesssim \int \varphi(f/\lambda) \, dx.
\end{align*}
The operator $\mathcal{M}_{\mathcal{S}}$ depends on the $(\varepsilon, d)$-configurations contained in $\mathcal{S}$,
and in particular it is always bounded (in $L^1$) for collections with a \emph{finite} number of rectangles.
Also, such weak-type inequalities do not hold for \emph{all} $(\varepsilon,d)$-configurations; some conditions of $\varepsilon$ and $d$ need to be specified.
For example, when $\varphi(x) = x^p$ then our results show that \eqref{questions:general_weak_type} holds uniformly
over all nonnegative functions $f$, all $\lambda \in (0, \infty)$, and all $(\varepsilon,d)$-configurations $\mathcal{S}$ with bounded $A_p(\mathcal{S})$.

\begin{question} \label{Q3}
  For which convex increasing functions $\varphi \colon [0,\infty) \to [0, \infty)$ satisfying $\varphi(0) = 0$
  does there exist an unbounded set $\Omega \subseteq (0,\frac{1}{2})\times \mathbb{N}$ such that
  \begin{align} \label{questions:general_weak_type}
    |\mathcal{M}_{\mathcal{S}}f > \lambda| \lesssim \int \varphi(f/\lambda) \, dx
  \end{align}
  holds uniformly over all nonnegative $f$, all $\lambda  \in (0, \infty)$, and all $(\varepsilon,d)$-configurations with $(\varepsilon,d) \in \Omega$?
\end{question}
Our results show that $\varphi(x) = x \log(e+x)$ is contained in the set of such functions, and that $\varphi(x) = x$ is not.
It would be interesting to know whether there exists a threshold function between $x$ and $x \log(e+x)$,
or whether one can formulate a condition that separates the two groups of functions from each other.

%\red{[L: summary from Darek?]}

%\bibliography{bibliography}

\begin{thebibliography}{KMPRR22}
\scriptsize
  \bibitem[BS88]{BS88} C.~Bennett, R.~Sharpley, 
 \emph{ Interpolation of operators}, Pure and Applied Mathematics 129, Academic Press, Inc., Boston, MA, 1988.


	\bibitem[CW71]{CW71} R.R.~Coifman, G.~Weiss, 
	\emph{Analyse harmonique non-commutative sur certains espaces homog{\`e}nes}, 
	Springer-Verlag, Berlin, 1971.

  \bibitem[C79]{C79} A.~C\'ordoba, \emph{Maximal functions, covering lemmas and Fourier multipliers,} Fourier
analysis in Euclidean spaces (Proc. Sympos. Pure Math., Williams
Coll., Williamstown, Mass. 1978), Part 1, 29--50, Proc. Sympos. Pure
Math. XXXV, Amer. Math. Soc., Providence, R.I., 1979.

    \bibitem[CF75]{CF75} A.~C\'ordoba, R.~Fefferman, 
	\emph{A geometric proof of the strong maximal theorem}, Ann. of Math.~{\bf 102} (1975), 95--100.
	
	
	\bibitem[Ed98]{Ed98} G.~Edgar, 
	\emph{Integral, Probability, and Fractal Measures},
	Springer-Verlag, New York, 1998.
	
%	\bibitem[Fe19]{Fe19} E.~Fern{\'a}ndez, 
%	\emph{An{\'a}lisis de Fourier en el toro infinito-dimensional},
%	Ph.D. thesis, Universidad de La Rioja, 2019.
	
%	\bibitem[FR19]{FR19} E.~Fern{\'a}ndez, L.~Roncal,  
%	\emph{On the absolute divergence of Fourier series on the infinite-dimensional torus},  Colloq. Math.~{\bf 157} (2019), 143--155.
	
	\bibitem[FR20]{FR20} E.~Fern{\'a}ndez, L.~Roncal,  
	\emph{A decomposition of Calder{\'o}n--Zygmund type and some observations on differentiation of integrals on the infinite-dimensional torus}, Potential Anal.~{\bf 53} (2020), 1449--1465.
	
%	\bibitem[Hy10]{Hy10} T.~Hyt\"onen, \emph{A framework for non-homogeneous analysis on metric spaces, and the RBMO space of Tolsa}, Publ. Mat.~{\bf 54} (2010), 485--504.

\bibitem[HJ05]{HJ05}
P.A.~Hagelstein, R.L.~Jones,
\emph{On restricted weak type $(1,1)$: the continuous case},
Proc. Amer. Math. Soc.~{\bf 133} (2005), 185--190
 
	\bibitem[Ko21]{Ko21}
	D.~Kosz, \emph{On differentiation of integrals in the infinite-dimensional torus}, Studia Math.~{\bf 258} (2021), 103–119.
	
	\bibitem[Ko22]{Ko22}
	D.~Kosz, \emph{On the doubling condition in the infinite-dimensional setting},  Bull. Aust. Math. Soc. (2023), published online.
	
	\bibitem[KMPRR22]{KMPRR22}
	D.~Kosz, J.~Mart{\'i}nez-Perales, V.~Paternostro, E.~Rela, L.~Roncal, \emph{Maximal operators on the infinite-dimensional torus}, Math. Ann. (2022), published online.
	
	\bibitem[La97]{La97}
	R.~Lata{\l}a, \emph{Estimation of moments of sums of independent real random variables}, Ann. Probab.~{\bf 25} (1997), 1502--1513.
	
	\bibitem[Mo05]{Mo05}
	 K.H.~Moon,
	 \emph{On restricted weak type $(1,1)$}, Proc. Amer. Math. Soc.~{\bf 4}2 (1974), 148--152.
	
	\bibitem[PS09]{PS09}
	M. Paluszy{\'n}ski, K.~Stempak, \emph{On quasi-metric and metric spaces}, Proc. Amer. Math. Soc.~{\bf 137} (2009), 4307--4312.
	
	\bibitem[RdF78]{RdF78}
	J.L.~Rubio~de~Francia, \emph{Nets of subgroups in locally compact groups}, Comment. Math. Prace Mat.~{\bf 20} (1977/78), 453--466.
	
	\bibitem[RdF80]{RdF80}
	J.L.~Rubio~de~Francia, \emph{Convergencia de series de Fourier de infinitas variables}, Publ. Sec. Mat. Univ. Aut{\`o}noma Barcelona~{\bf 21} (1980), 237--241.
	
	\bibitem[Zy67]{Zy67}
	A.~Zygmund, \emph{A note on the differentiability of integrals}, Colloq. Math.~{\bf 16} (1967), 199--204.
		
\end{thebibliography}
%\bibliographystyle{abbrv}

\subsection*{Acknowledgment} The authors would like to thank the referee for pointing out useful literature related to the topic, as well as for several remarks that improved the presentation of the results.   

\subsection*{Funding}
D.K. and L.R. were supported by the Basque Government (BERC 2022-2025) and by the Spanish State Research Agency (CEX2021-001142-S).
D.K. was also supported by the Foundation for Polish Science (START 032.2022), and L.R. was also supported by the Spanish State Research Agency (PID2020-113156GB-I00/AEI/10.13039/501100011033 with acronym ``HAPDE'' and  RYC2018-025477-I), and by Ikerbasque.
\\
G.R. was supported by PID2019-105599GB-I00/AEI/10.13039/501100011033.

\end{document}